\documentclass[sn-mathphys,Numbered]{sn-jnl}

\usepackage{graphicx}
\usepackage{multirow}
\usepackage{amsmath,amssymb,amsfonts}
\usepackage{amsthm}
\usepackage{mathrsfs}
\usepackage[title]{appendix}
\usepackage{xcolor}
\usepackage{textcomp}
\usepackage{manyfoot}
\usepackage{booktabs}
\usepackage{algorithm}
\usepackage{algorithmicx}
\usepackage{algpseudocode}
\usepackage{listings}
\usepackage[latin1]{inputenc}
\usepackage{makecell}
\usepackage{tabularray}
\usepackage{tikz-cd}
\usepackage{mathbbol}
\usepackage{float}
\usepackage{caption}
\usepackage{array}

\theoremstyle{thmstyleone}
\newtheorem{theorem}{Theorem}[section]

\theoremstyle{thmstyletwo}

\newtheorem{conjecture}[theorem]{Conjecture}

\theoremstyle{thmstylethree}

\everymath{\displaystyle}

\begin{document}

\title[On Thom Polynomials for $\mathcal{A}$-singularities of parametric curves and degeneracy loci formulas]{On Thom Polynomials for $\mathcal{A}$-singularities of parametric curves and degeneracy loci formulas}

\author*[1]{\fnm{João Vítor} \sur{Pissolato}}\email{jvpissolato@usp.br}

\affil*[1]{\orgdiv{Department of Mathematics}, \orgname{Instituto de Ciências Matemáticas e de Computação}, \orgaddress{\street{Trabalhador São-Carlense Avenue, 400}, \city{São Carlos}, \postcode{13566-590}, \state{São Paulo}, \country{Brazil}}}

\abstract{In Singularity Theory of smooth mappings, an important direction is the local study of singular points, commonly known as germs of mappings. One topological object that carries geometrical information about those singularities is the Thom Polynomial, which is a universal expression in Chern classes that depends only on the singularity type.
	
The Thom Polynomial has been studied for $\mathcal{K}$-singularities by many authors, and diverse applications, such as obtaining enumerative formulas for $0$-stable invariants, have been discovered, making it an interesting object to study further. For $\mathcal{K}$-singularities, the Thom Polynomial is written in terms of quotient Chern classes, as long as the relative codimension is fixed. On the other hand, for $\mathcal{A}$-singularities, it is simply written in terms of the Chern classes of the source and those of the target, and, in many cases, it cannot be expressed through quotient Chern classes.
	
In this research article, we make an in-depth investigation of the computations for Thom Polynomials for $\mathcal{A}$-singularities of parametric curves germs in the plane and space cases using the framework established by Ohmoto and the Restriction Method introduced by Rim\'{a}nyi, and we also propose a geometric way for such computations using degeneracy loci Chern classes formulas proven by Fulton, which seems a promising way to generalize the Thom Polynomial expression for the stabilization of such singularities in $\mathbb{C}^N$, depending exclusively on $N$ as a parameter.}

\keywords{Curves singularities, Thom Polynomials, Degeneracy loci}

\maketitle	

\section{Introduction}\label{sec1}

In Singularity Theory, one way to study singular points of a mapping is to study its orbit under a action that comes from a group in the space of germs. The most usual ones are the Mather groups, mainly the $\mathcal{A}$ and $\mathcal{K}$ groups. In a certain sense, the $\mathcal{A}$ group induces an action that changes coordinates in the source and in the target, and this equivalence relation is called ``right-left equivalence''. On the other hand, the $\mathcal{K}$ group induces an equivalence that is called ``contact equivalence'', and it detects the different contacts that the germs can have at the origin.

Given a stable $\mathcal{K}$-type singularity $\xi$ in $\mathcal{O}_0(m,n)$, define the relative codimension $\kappa=n-m$, and given any mapping $f:M\to N$ with relative codimension $\kappa$, one can consider the set that consists of the points in $M$ such that $f$ has a singularity of type $\eta$ at $x$, denoted by $\eta(f)$. Then, considering the fundamental class of its closure in the Borel-Moore homology, $[\overline{\eta(f)}]$, and applying the Poincaré Duality, one gets a cohomology class in the singular cohomology. In fact, this cohomology class is written as a polynomial expression of quotient Chern classes of $f$ with integers coefficients, and this polynomial is a universal expression that does not depend on $f$. This is the Thom Polynomial for the $\mathcal{K}$-type singularity $\eta$. We refer the works \cite{11}, \cite{18} and \cite{1} for more details.

If we change the context to $\mathcal{A}$-singularities, more singularity types appear, since the action group is smaller. In this case, a Thom Polynomial theory also exists, which has been shown originally by Haefliger and Kosinski in \cite{21}, and it is simply written in terms of Chern classes of the source and those of the target. More details about this context can be found in \cite{18}.

One difficult question is to determine the explicit expression of the Thom Polynomial for a given singularity type. Over the years, many authors have obtained some Thom Polynomials. A classical way is to find a suitable resolution of the singular locus, however it is a very difficult task. More recently, in \cite{3}, Rim\'{a}nyi introduced an effective way to compute Thom Polynomials, especially for singularities of small codimension and with quasi-homogeneous normal form, called the ``Restriction Method''. This technique makes use of classification tables of singularities and the existence theorem of the Thom Polynomial. A reasonable framework for Thom Polynomial for $\mathcal{A}$-singularities was proposed by Ohmoto and Sasajima in \cite{2}, which makes use of the Restriction Method to compute the Thom Polynomials of unstable $\mathcal{A}$-singularity types.

In \cite{15}, Kazarian, Lando and Zvonkine obtained a universal cohomology expressions for singularities of family of mappings between curves, that could be used to study the $(\mathbb{C},0)\to(\mathbb{C},0)$ case. In the present paper, we give Thom Polynomials for $\mathcal{A}$-type singularities of plane and space curves germs $(\mathbb{C},0)\to(\mathbb{C}^n,0)$, for $n=2,3$, up to some codimension, where we start to face some technical problems which makes it difficult to proceed further. Also, we propose some expressions for an arbitrary $n$ using the degeneracy loci Chern classes formulas obtained by Fulton in \cite{9}. Furthermore, this method seems to have some connection with the techniques seen in the works of Buch, Feh\'{e}r, Fulton and Rim\'{a}nyi in \cite{16}, \cite{22} and \cite{17} for the Thom Polynomial of quiver varieties, where they obtained a formula expressing the cohomology class as a linear combination of products of Schur determinants. When interpreted for Thom Polynomials, this formula has the form $$Tp_r=\sum_\lambda c_\lambda(r)s_{\lambda_1}(x^1:x^0)s_{\lambda_2}(x^2:x^1)\cdots s_{\lambda_n}(x^n:x^{n-1}),$$ where $c_\lambda(r)$ is a certain integer and each $s_{\lambda_i}$ is a Schur determinant.

The author is hopeful that such method for the Thom Polynomial of $\mathcal{A}$-singularities of parametric curves can be generalized, giving a new direction to the Thom Polynomial theory and their computations.

\section{Basic notions and results about $\mathcal{A}$-singularities}\label{sec2}

We denote by $\mathcal{O}_0(m,n)$ the set of all germs of smooth mappings $f:(\mathbb{C}^m,0)\to(\mathbb{C}^n,0)$, and by $\mathcal{A}$ the set of all pairs $(h,k)$, where $h:(\mathbb{C}^m,0)\to(\mathbb{C}^m,0)$ and $k:(\mathbb{C}^n,0)\to(\mathbb{C}^n,0)$ are germs of smooth diffeomorphisms. It is not difficult to see that $\mathcal{A}$ is a group under composition and it acts on $\mathcal{O}_0(m,n)$ via $$(h,k)\cdot f=k\circ f \circ h^{-1},$$ which can be interpreted as a smooth change of coordinates on source and target.

Classically, the characteristic classes of a triangulated smooth manifold have been defined as the obstruction for extending a smooth vector field without singularities that is already defined on a certain stratum to a higher stratum. There are other definitions and they are all equivalent. Here, we make use of the Chern classes for a vector bundle, which are sufficient for our computations. We refer the reader to see \cite{10} for more details.

Given two vector bundles $E$ and $F$, the quotient total Chern class is defined as $$c(E-F)=\sum_{i\geq0}c_i(E-F)=\frac{1+c_1(E)+c_2(E)+\cdots}{1+c_1(F)+c_2(F)+\cdots}.$$

However, to our interests, we need the expression of the top Chern class. If \linebreak$p:E\to M$ is a complex vector bundle of rank $r$, the top Chern class of $E$ is defined as $$c_r(E)=(p^\ast)^{-1}(\text{Dual}([Z]))\in H^{2r}(M;\mathbb{Z}),$$ where $\text{Dual}$ is the Poincaré Duality, that is, the inverse of the isomorphism $$\varphi\mapsto\varphi\frown[M],$$ and $Z$ is the zero section of $E$.

Given a stable $\mathcal{K}$-type singularity $\eta$ in $\mathcal{O}_0(m,n)$, define the relative codimension as $\kappa=n-m$. The Thom Polynomial of singularity type $\eta$ is the unique universal expression in Chern classes $Tp(\eta)\in\mathbb{Z}[c_1, c_2, \ldots]$ such that, for any given smooth mapping $f:M\to N$ with relative codimension $\kappa$, denoting by $\eta(f)$ the set of points $x\in M$ such that $f$ has a singularity of type $\eta$ at $x$, we have $$\text{Dual}([\overline{\eta(f)}])=Tp(\eta)(c(f))\in H^{2\text{codim}(\eta)}(M),$$ where $c_i(f)=c_i(f^\ast TN-TM)$.

The framework for $\mathcal{A}$-singularities, as described in \cite{2}, is the following: Let $\eta$ be an unstable $\mathcal{A}$-singularity type in $\mathcal{O}_0(m,n)$. Given two vector bundles $E$ and $F$ over a holomorphic manifold $M$, with rank $m$ and $n$, respectively, and a generic family of map-germs $\varphi:(E,Z_E)\to(F,Z_F)$ such that the diagram $$\begin{tikzcd}
	(E,Z_E)\arrow{rr}{\varphi}\arrow{rd}[swap]{\pi_E}&&(F,Z_F)\arrow{ld}{\pi_F}\\&M&
\end{tikzcd}$$ commutes, we associate the cohomology class dual to the fundamental class of the closure of the set $\eta(\varphi)$. Then, for an appropriately generic $\varphi$, that cohomology class is expressed by $$\text{Dual}([\overline{\eta(\varphi)}])=Tp_\mathcal{A}(\eta)(a_1,\ldots,a_m,b_1,\ldots,b_n)\in H^{2\text{codim}(\eta)}(M),$$ where $a_i=c_i(E)$ and $b_j=c_j(F)$. In general, such Thom Polynomial is written in terms of the Chern classes of the source and those of the target. However, if the $\mathcal{A}$-orbit of $\eta$ is dense in its $\mathcal{K}$-orbit, then the above Thom Polynomial can be rewritten in terms of quotient Chern classes.

\section{Thom Polynomial computations}\label{sec3}

In this section, we present our computational results of the Thom Polynomial for some $\mathcal{A}$-type singularities of parametric curves in the plane and in the space, based on the classification made by many authors. The computations shown in the present section can be used to justify our results in Section \ref{sec4}, where we give an alternative geometric interpretation for those computations.

\subsection{Thom Polynomials for parametric plane curves}

In \cite{4}, Bruce and Gaffney obtained the classification for simple singularities of complex parametric plane curves. A part of their classification is given in Table \ref{tab1} and we use the dash notation to distinguish the singularities in the same family given by Bruce and Gaffney, but with non quasi-homogeneous normal forms, and we also give their $\mathcal{A}$-codimension and $\mathcal{A}_e$-miniversal deformation, up to codimension $7$.

{\renewcommand{\arraystretch}{2}\begin{table}
		\centering
		\caption{Parametric plane curve singularities of $\mathcal{A}$-cod $\leq7$}
		\label{tab1}
		\begin{tabular}{|c|c|c|c|}\hline
			\text{Symbol} & \text{Normal form} & $\mathcal{A}$-cod & $\mathcal{A}_e$-\text{miniversal deformation}\\\hline
			$A_2$ & $(t^2,t^3)$ & $2$ & $(t^2,t^3+u_1t)$\\\hline
			$A_4$ & $(t^2,t^5)$ & $3$ & $(t^2,t^5+u_1t+u_2t^3)$\\\hline
			$A_6$ & $(t^2,t^7)$ & $4$ & $(t^2,t^7+u_1t+u_2t^3+u_3t^5)$\\\hline
			$E_6$ & $(t^3,t^4)$ & $4$ & $(t^3+u_1t,t^4+u_2t+u_3t^2)$\\\hline
			$A_8$ & $(t^2,t^9)$ & $5$ & $(t^2,t^9+u_1t+u_2t^3+u_3t^5+u_4t^7)$\\\hline
			$E_8$ & $(t^3,t^5)$ & $5$ & $(t^3+u_1t,t^5+u_2t+u_3t^2+u_4t^4)$\\\hline
			$A_{10}$ & $(t^2,t^{11})$ & $6$ & $(t^2,t^{11}+u_1t+u_2t^3+u_3t^5+u_4t^7+u_5t^9)$\\\hline
			$E_{12}'$ & $(t^3,t^7+t^8)$ & $6$ & $(t^3+u_1t,t^7+t^8+u_2t+u_3t^2+u_4t^4+u_5t^5)$\\\hline
			$W_{12}'$ & $(t^4,t^5+t^7)$ & $6$ & $(t^4+u_1t+u_2t^2,t^5+t^7+u_3t+u_4t^2+u_5t^3)$\\\hline
			$A_{12}$ & $(t^2,t^{13})$ & $7$ & $(t^2,t^{13}+u_1t+u_2t^3+u_3t^5+u_4t^7+u_5t^9+u_6t^{11})$\\\hline
			$E_{12}$ & $(t^3,t^7)$ & $7$ & $(t^3+u_1t,t^7+u_2t+u_3t^2+u_4t^4+u_5t^5+u_6t^8)$\\\hline
			$E_{14}'$ & $(t^3,t^8+t^{10})$ & $7$ & $(t^3+u_1t,t^8+t^{10}+u_2t+u_3t^2+u_4t^4+u_5t^5+u_6t^7)$\\\hline
			$W_{12}$ & $(t^4,t^5)$ & $7$ & $(t^4+u_1t+u_2t^2,t^5+u_3t+u_4t^2+u_5t^3+u_6t^7)$\\\hline
			$W_{1,1}^\#$ & $(t^4,t^6+t^7)$ & $7$ & \makecell{$(t^4+u_1t+u_2t^2,t^6+t^7+u_3t$\\$+u_4t^2+u_5t^3+u_6t^5)$}\\\hline
		\end{tabular}
\end{table}}

We use the common notation in Thom-Mather theory and more details can be seen in \cite{20}. Given an integer $k\geq 1$, we denote by $\theta_k$ the $\mathcal{O}_{\mathbb{C}^k,0}$-module of germs of vector fields at the origin in $\mathbb{C}^k$. For a germ $f:(\mathbb{C},0)\to(\mathbb{C}^n,0)$, let $\theta_f$ be the $\mathcal{O}_{\mathbb{C}^n,0}$-module of vector fields along $f$, and $\mathfrak{m}_n$ be the maximal ideal of the local ring $\mathcal{O}_{\mathbb{C}^n,0}$. To compute the $\mathcal{A}$-codimension of $f$, we simply apply the definition $$\mathcal{A}-\textnormal{codim}(f)=\textnormal{dim}_\mathbb{C}\left(\frac{\mathfrak{m}_1\theta_f}{T\mathcal{A}\cdot f}\right),$$ where $T\mathcal{A}\cdot f$ is the tangent space to the orbit $\mathcal{A}\cdot f$ defined as $$T\mathcal{A}\cdot f=tf(\mathfrak{m}_1\theta_1)+\omega f(\mathfrak{m}_n\theta_n),$$ where $tf$ and $\omega f$ are the mappings defined by $tf(\xi)=df\circ\xi$ and $\omega f(\eta)=\eta\circ f$, for $\xi\in\theta_1$ and $\eta\in\theta_n$. For example, for the $E_8$ singularity, we have $$T\mathcal{A}\cdot E_8=\left\langle\left(\begin{array}{c}3t^3\\5t^5\end{array}\right)\right\rangle+\left\{\left.\left(\begin{array}{c}\eta_1(t^3,t^5)\\\eta_2(t^3,t^5)\end{array}\right)\right|\eta_1,\eta_2\in\mathfrak{m}_2\theta_2\right\},$$ and we must pay attention to the mixed module structure in it. It is not difficult to see that the only missing monomials are $$\left(\begin{array}{c}t\\0\end{array}\right),\left(\begin{array}{c}t^2\\0\end{array}\right),\left(\begin{array}{c}0\\t\end{array}\right),\left(\begin{array}{c}0\\t^2\end{array}\right),\left(\begin{array}{c}0\\t^4\end{array}\right)$$ and that they are linearly independent in $\frac{\mathfrak{m}_1\theta_{E_8}}{T\mathcal{A}\cdot E_8}$. Therefore, it has $\mathcal{A}$-codimension $5$.

For the $\mathcal{A}_e$-miniversal deformation, it is a similar computation with the extended tangent space to the $\mathcal{A}$-orbit $$T\mathcal{A}_e\cdot E_8=tf(\theta_1)+\omega f(\theta_2)=\left\langle\left(\begin{array}{c}3t^2\\5t^4\end{array}\right)\right\rangle+\left\{\left.\left(\begin{array}{c}\eta_1(t^3,t^5)\\\eta_2(t^3,t^5)\end{array}\right)\right|\eta_1,\eta_2\in\theta_2\right\}.$$ The missing monomials are the same as above, but this time, the second and fifth monomials are linearly dependent in the respective quotient of $\theta_{E_8}$ by ${T\mathcal{A}_e\cdot E_8}$, so we can discard the second one, obtaining the $\mathcal{A}_e$-miniversal deformation shown in Table \ref{tab1}.

Let us compute the Thom Polynomial for the $E_8$-singularity by using the Restriction Method introduced by Rimányi in \cite{3}. Since $\mathcal{A}-\textnormal{cod}(E_8)=5$, we write a general homogeneous polynomial expression of degree $5$ in Chern classes of the source and those of the target \begin{equation}\label{exp2}\begin{array}{rcl}Tp_\mathcal{A}(E_8)&=&x_1a_1^5+x_2a_1^4b_1+x_3a_1^3b_1^2+x_4a_1^3b_2+x_5a_1^2b_1^3+x_6a_1^2b_1b_2+x_7a_1b_1^4\\&&+x_8a_1b_1^2b_2+x_9a_1b_2^2+x_{10}b_1^5+x_{11}b_1^3b_2+x_{12}b_1b_2^2,\end{array}\end{equation} and we need to determine the unknown integers $x_i$, for $1\leq i\leq 12$.

First, we get the equation that Rim\'{a}nyi calls the ``principal equation''. Since the $\mathcal{A}_e$-miniversal deformation of $E_8$ is $$\widetilde{E_8}(t,u_1,u_2,u_3,u_4)=(t^3+u_1t,t^5+u_2t+u_3t^2+u_4t^4),$$ which is quasi-homogeneous, we consider the algebraic torus $T=\mathbb{C}^\ast$ and define actions in the source and in the target that keep $\widetilde{E_8}$ invariant. For example, consider $$\rho_1=\alpha\oplus\alpha^2\oplus\alpha^4\oplus\alpha^3\oplus\alpha,\text{ }\rho_2=\alpha^3\oplus\alpha^5,\text{ for }\alpha\in T,$$ which satisfies $\widetilde{E_8}\circ\rho_1=\rho_2\circ \widetilde{E_8}$. Consider the dual tautological line bundle $\ell=\mathcal{O}_{\mathbb{P}_\mathbb{C}^N}(1)$ over some complex projective space of large enough dimension $\mathbb{P}_\mathbb{C}^N$ and the vector bundles according to the previously defined actions $$E_1=\ell,\text{ }E_1'=\ell^{\otimes2}\oplus\ell^{\otimes4}\oplus\ell^{\otimes3}\oplus\ell,\text{ }E_2=\ell^{\otimes3}\oplus\ell^{\otimes5}.$$ We have that $$H^\ast(\mathbb{P}_\mathbb{C}^N)\cong\frac{\mathbb{Z}[a]}{\langle a^{N+1}\rangle},$$ where $a$ is the class represented by a generic hyperplane in the projective space. Then we compute the total Chern classes \begin{align*}
	c(E_1)&=1+a\\
	c(E_1\oplus E_1')&=(1+a)(1+2a)(1+4a)(1+3a)(1+a)\\
	&=1+11a+45a^2+85a^3+74a^4+24a^5\\
	c(E_2)&=(1+3a)(1+5a)\\
	&=1+8a+15a^2.
\end{align*}

Now, trivializing the vector bundles and gluing the normal form of $E_8$ along the trivializing neighborhoods, we get a family of holomorphic map-germs $f_{E_8}$
$$\begin{tikzcd}
	E_1\oplus E_1'\arrow{rr}{f_{E_8}}\arrow{rd}&&E_2\arrow{ld}\\&\mathbb{P}_\mathbb{C}^N&
\end{tikzcd}$$ whose restriction to each fiber is $\mathcal{A}$-equivalent to $E_8$, called the ``universal map for $E_8$''. Then, its corresponding Chern classes are $$a_1(f_{E_8})=a,\text{ }b_1(f_{E_8})=8a\text{ and }b_2(f_{E_8})=15a^2.$$ We plug in these Chern classes in the right-hand side of \ref{exp2} to get one side of the principal equation. On the other hand, we have that $E_8$ has only the origin as $E_8$ singularity, so $\overline{E_8(f_{E_8})}=\text{Z}(E_1\oplus E_1')$. Then $$Tp_\mathcal{A}(f_{E_8})=\text{Dual}([\overline{E_8(f_{E_8})}])=\text{Dual}([\text{Z}(E_1\oplus E_1')])=c_5(E_1\oplus E_1')=24a^5$$ and we obtain from \ref{exp2} the equation \begin{align*}x_1 + 8x_2 + 64x_3 + 15x_4 + 512x_5 + 120x_6 + 4096x_7 + 960x_8 + 225x_9 &\\ + 32768x_{10} + 7680x_{11} + 1800x_{12}&=24\end{align*} The same principal equation can be obtained if we choose the other $\mathcal{A}_e$-miniversal deformation representative, since the powers of the algebraic torus actions will be the same, up to a permutation of their ordering, so the total Chern classes involved will be the same.

Now we apply the same method to other $\mathcal{A}$-singularities of $\mathcal{A}$-codimension less than or equal to $5$ to get the equations that Rim\'{a}nyi calls the ``homogeneous equations''. This name comes from the fact that, as shown in \cite{3}, they are all equal to $0$, since the $E_8$ singularity does not appear in their $\mathcal{A}_e$-miniversal deformations, including the immersion germ $A_0(t)=(t,0)$. Actually, we can get six independent equations from $A_0$, because of the zero second coordinate, and one independent equation from each singularity $A_2$, $A_4$, $A_6$, $A_8$ and $E_6$. One obtainable linear system is given in Table \ref{tab2}.

\begin{table}
	\centering
	\caption{A system of linear equations that determine $Tp_\mathcal{A}(E_8)$}
	\label{tab2}
	{\renewcommand{\arraystretch}{2}\begin{tabular}{|ccc|}\hline
	\makecell{$x_1 + 8x_2 + 64x_3 + 15x_4 + 512x_5 + 120x_6 + 4096x_7$\\$+ 960x_8 + 225x_9 + 32768x_{10} + 7680x_{11} + 1800x_{12}$}&=&24\\\hline
	$x_{10}$&=&0\\\hline
	$x_1  + x_2 + x_3 + x_5 + x_7+ x_{10}$&=&0\\\hline
	\makecell{$x_1 + 2x_2 + 4x_3 + x_4 + 8x_5 + 2x_6 + 16x_7 + 4x_8 + x_9+ 32x_{10}+ 8x_{11} + 2x_{12}$}&=&0\\\hline
	\makecell{$32x_1 + 48x_2 + 72x_3 + 16x_4 + 108x_5 + 24x_6$\\$+ 162x_7 + 36x_8 + 8x_9 + 243x_{10} + 54x_{11} + 12x_{12}$}&=&0\\\hline
	\makecell{$x_1 + 3x_2 + 9x_3 + 2x_4 + 27x_5 + 6x_6 + 81x_7$\\$+ 18x_8 + 4x_9 + 243x_{10} + 54x_{11} + 12x_{12}$}&=&0\\\hline
	\makecell{$x_1 + 4x_2 + 16x_3 + 3x_4 + 64x_5 + 12x_6 + 256x_7$\\$+ 48x_8 + 9x_9 + 1024x_{10} + 192x_{11} + 36x_{12}$}&=&0\\\hline
	\makecell{$x_1 + 5x_2 + 25x_3 + 6x_4 + 125x_5 + 30x_6 + 625x_7$\\$+ 150x_8 + 36x_9 + 3125x_{10} + 750x_{11} + 180x_{12}$}&=&0\\\hline
	\makecell{$x_1 + 7x_2 + 49x_3 + 10x_4 + 343x_5 + 70x_6 + 2401x_7$\\$+ 490x_8 + 100x_9 + 16807x_{10} + 3430x_{11} + 700x_{12}$}&=&0\\\hline
	\makecell{$x_1 + 9x_2 + 81x_3 + 14x_4 + 729x_5 + 126x_6 + 6561x_7$\\$+ 1134x_8 + 196x_9 + 59049x_{10} + 10206x_{11} + 1764x_{12}$}&=&0\\\hline
	\makecell{$x_1 + 11x_2 + 121x_3 + 18x_4 + 1331x_5 + 198x_6 + 14641x_7$\\$+ 2178x_8 + 324x_9 + 161051x_{10} + 23958x_{11} + 3564x_{12}$}&=&0\\\hline
	\makecell{$x_1 + 7x_2 + 49x_3 + 12x_4 + 343x_5 + 84x_6 + 2401x_7$\\$+ 588x_8 + 144x_9 + 16807x_{10} + 4116x_{11} + 1008x_{12}$}&=&0\\\hline
\end{tabular}}
\end{table} By solving it, we obtain $$Tp_\mathcal{A}(E_8)=-28a_1^5 + 46a_1^4b_1 - 20a_1^3b_1^2 - 35a_1^3b_2 + 2a_1^2b_1^3 + 26a_1^2b_1b_2 - 3a_1b_1^2b_2 - 7a_1b_2^2 + b_1b_2^2.$$

In a similar way, all the Thom Polynomial for singularities up to $\mathcal{A}$-codimension 5 can be computed. At $\mathcal{A}$-codimension $6$, singularities with non quasi-homogeneous normal forms appears. To compute its Thom Polynomial, we restrict it to a suitable jet space. For example, to compute $Tp_\mathcal{A}(E_{12}')$, we consider the restriction of its $\mathcal{A}_e$-miniversal deformation to the $7$-jet space, getting the germ $$(t^3+u_1t,t^7+u_2t+u_3t^2+u_4t^4+u_5t^5),$$ which is quasi-homogeneous and we can get the principal equation.  We use the same argument to get the homogeneous equation that comes from $W_{12}'$, which we restrict to the $5$-jet space.

It is important to note that, with the Restriction Method and the above remark, we can compute the Thom Polynomials up to codimension $6$. However, at codimension $7$, the singularity $W_{1,1}^\#$ appears, and we would need to get an restriction equation from it. Its normal form is not quasi-homogeneous, and if we restrict it to the $6$-jet space, we get the germ $(t^4,t^6)$, which has infinite $\mathcal{A}$-codimension, so it is not possible to apply the Thom Polynomial theorem and we cannot get an equation from that singularity. Therefore, to obtain the Thom Polynomial for the $\mathcal{A}$-codimension $7$ singularities, the last homogeneous equation comes from $$\widetilde{W}_{18}(t,u)=(t^4+u_1t+u_2t^2,t^7+u_3t+u_4t^2+u_5t^3+u_6t^5+u_7t^6+u_8t^9+u_9t^{13}),$$ which has no singularities of type $A_{12}$, $E_{12}$, $E_{14}'$ or $W_{12}$.

Carrying out those computations, we obtain the Thom Polynomial expressions for complex plane curves $\mathcal{A}$-singularities. Thanks to a remark from  Ohmoto, we noticed that such expressions can be factorized. Therefore, we conclude the following result, where we give the Thom Polynomials in a factorized form.

\begin{theorem}\label{teo1}
	
The Thom Polynomials for $\mathcal{A}$-singularities of germs $(\mathbb{C},0)\to(\mathbb{C}^2,0)$ of $\mathcal{A}$-codimension less than or equal to $7$, except for $W_{1,1}^\#$, are given in Table \ref{tab3}.

\begin{table}
	\centering
	\caption{Factorized Thom Polynomials for $\mathcal{A}\text{-cod}\leq7$ parametric plane curve singularities, except for $W_{1,1}^\#$}
	\label{tab3}
	{\renewcommand{\arraystretch}{2}\begin{tabular}{|c|c|}\hline
	\text{Symbol} & \text{Thom Polynomial}\\\hline
		$A_2$ & $c_2=a_1^2 - a_1b_1 + b_2$\\\hline
		$A_4$ & $c_2(-5a_1 + b_1)$\\\hline
		$A_6$ & $c_2(-5a_1 + b_1)(-7a_1 + b_1)$\\\hline
		$E_6$ & \makecell{$c_1c_3+c_2^2+2c_4$\\$=c_2(4a_1^2 - 2a_1b_1 + b_2)$}\\\hline
		$A_8$ & $c_2(-7a_1 + b_1)(49a_1^2 - 16a_1b_1 + b_1^2 + b_2)$\\\hline
		$E_8$ & $(c_1c_3+c_2^2+2c_4)(-7a_1 + b_1)$\\\hline
		$A_{10}$ & \makecell{$c_2(4033a_1^4 - 2316a_1^3b_1 + 454a_1^2b_1^2 - 36a_1b_1^3 + b_1^4$\\$+ 146a_1^2b_2 - 38a_1b_1b_2 + 2b_1^2b_2 + b_2^2)$}\\\hline
		$E_{12}'$ & $(c_1c_3+c_2^2+2c_4)(47a_1^2 - 12a_1b_1 + b_1^2 - b_2)$\\\hline
		$W_{12}'$ & \makecell{$2c_1^2c_4 + 3c_1c_2c_3 + 10c_1c_5 + c_2^3+ 7c_2c_4 + c_3^2 + 12c_6$ \\$=(c_1c_3+c_2^2+2c_4)(9a_1^2 - 3a_1b_1 + b_2)$}\\\hline
		$A_{12}$ & \makecell{$c_2(-11a_1 + b_1)(5167a_1^4 - 2982a_1^3b_1 + 568a_1^2b_1^2$\\$- 42a_1b_1^3 + b_1^4 + 272a_1^2b_2 - 70a_1b_1b_2 + 4b_1^2b_2 + b_2^2)$}\\\hline
		$E_{12}$ & $(c_1c_3+c_2^2+2c_4)(47a_1^2 - 12a_1b_1 + b_1^2 - b_2)(-11a_1 + b_1)$\\\hline
		$E_{14}'$ & \makecell{$(c_1c_3+c_2^2+2c_4)(-389a_1^3 + 131a_1^2b_1 - 19a_1b_1^2$\\$+ b_1^3 + 19a_1b_2 - 2b_1b_2)$}\\\hline
		$W_{12}$ & $(2c_1^2c_4 + 3c_1c_2c_3 + 10c_1c_5 + c_2^3+ 7c_2c_4 + c_3^2+ 12c_6)(-11a_1 + b_1)$\\\hline
	\end{tabular}}
\end{table}
\end{theorem}

\subsection{Thom Polynomials for parametric space curves}

In this section, we make the same computations as in the previous section for parametric space curves singularities. In \cite{6}, Gibson and Hobbs obtained the classification of simple singularities of parametric space curves. Naturally, all the singularities of the plane case appears with a $0$ in the third coordinate, which the authors call ``planar curves'', and other additional normal forms, which are called ``spatial curves''. Therefore, the same technical difficulty regarding the singularity $(t^4,t^6+t^7)$ is expected to occur in this setting. However, the Thom Polynomial computation becomes problematic earlier than that, at $\mathcal{A}$-codimension $12$. The singularities that we worked with are given in Table \ref{tab4}.

\begin{table}
	\centering
	\caption{Parametric space curves singularities of $\mathcal{A}\text{-cod}\leq12$}
	\label{tab4}
	{\renewcommand{\arraystretch}{2}\begin{tabular}{|c|c|c|}\hline
		\text{Normal form} & $\mathcal{A}$-cod & $\mathcal{A}_e$-\text{miniversal deformation}\\\hline
		$(t^2,t^3,0)$ & $3$ & $(t^2,t^3+u_1t,u_2t)$\\\hline
		$(t^2,t^5,0)$ & $5$ & $(t^2,t^5+u_1t+u_2t^3,u_3t+u_4t^3)$\\\hline
		$(t^3,t^4,t^5)$ & $6$ & $(t^3+u_1t,t^4+u_2t+u_3t^2,t^5+u_4t+u_5t^2)$\\\hline
		$(t^2,t^7,0)$ & $7$ & $(t^2,t^7+u_1t+u_2t^3+u_3t^5,u_4t+u_5t^3+u_6t^5)$\\\hline
		$(t^3,t^4,0)$ & $7$ & $(t^3+u_1t,t^4+u_2t+u_3t^2,u_4t+u_5t^2+u_6t^5)$\\\hline
		$(t^3,t^5,t^7)$ & $8$ & \makecell{$(t^3+u_1t,t^5+u_2t+u_3t^2+u_4t^4,$\\$t^7+u_5t+u_6t^2+u_7t^4)$}\\\hline
		$(t^2,t^9,0)$ & $9$ & \makecell{$(t^2,t^9+u_1t+u_2t^3+u_3t^5+u_4t^7,$\\$u_5t+u_6t^3+u_7t^5+u_8t^7)$}\\\hline
		$(t^3,t^5,0)$ & $9$ & \makecell{$(t^3+u_1t,t^5+u_2t+u_3t^2+u_4t^4,$\\$u_5t+u_6t^2+u_7t^4+u_8t^7)$}\\\hline
		$(t^4,t^5,t^6)$ & $9$ & \makecell{$(t^4+u_1t+u_2t^2,t^5+u_3t+u_4t^2+u_5t^3,$\\$t^6+u_6t+u_7t^2+u_8t^3)$}\\\hline
		$(t^3,t^7,t^8)$ & $10$ & \makecell{$(t^3+u_1t,t^7+u_2t+u_3t^2+u_4t^4+u_5t^5,$\\$t^8+u_6t+u_7t^2+u_8t^4+u_9t^5)$}\\\hline
		$(t^4,t^5,t^7)$ & $10$ & \makecell{$(t^4+u_1t+u_2t^2,t^5+u_3t+u_4t^2+u_5t^3,$\\$t^7+u_6t+u_7t^2+u_8t^3+u_9t^6)$}\\\hline
		$(t^2,t^{11},0)$ & $11$ & \makecell{$(t^2,t^{11}+u_1t+u_2t^3+u_3t^5+u_4t^7+u_5t^9,$\\$u_6t+u_7t^3+u_8t^5+u_9t^7+u_{10}t^9)$}\\\hline
		$(t^3,t^7+t^8,t^{11})$ & $11$ & \makecell{$(t^3+u_1t,t^7+t^8+u_2t+u_3t^2+u_4t^4+u_5t^5,$\\$t^{11}+u_6t+u_7t^2+u_8t^4+u_9t^5+u_{10}t^7)$}\\\hline
		$(t^4,t^5+t^7,t^{11})$ & $11$ & \makecell{$(t^4+u_1t+u_2t^2,t^5+t^7+u_3t+u_4t^2+u_5t^3,$\\$t^{11}+u_6t+u_7t^2+u_8t^3+u_9t^5+u_{10}t^6)$}\\\hline
		$(t^4,t^6,t^7)$ & $ 11$ & \makecell{$(t^4+u_1t+u_2t^2, t^6+u_3t+u_4t^2+u_5t^3+u_6t^5,$\\$t^7+u_7t+u_8t^2+u_9t^3+u_{10}t^5)$}\\\hline
		$(t^3,t^7+t^8,0)$ & $12$ & \makecell{$(t^3+u_1t,t^7+t^8+u_2t+u_3t^2+u_4t^4+u_5t^5,$\\$u_6t+u_7t^2+u_8t^4+u_9t^5+u_{10}t^7+u_{11}t^{10})$} \\\hline
		$(t^4,t^5+t^7,0)$ & $12$ & \makecell{$(t^4+u_1t+u_2t^2,t^5+t^7+u_3t+u_4t^2+u_5t^3,$\\$u_6t+u_7t^2+u_8t^3+u_9t^5+u_{10}t^6+u_{11}t^9)$}\\\hline
		$(t^3,t^7,t^{11})$ & $12$ & \makecell{$(t^3+u_1t,t^7+u_2t+u_3t^2+u_4t^4+u_5t^5+u_6t^8,$\\$t^{11}+u_7t+u_8t^2+u_9t^4+u_{10}t^5+u_{11}t^8)$}\\\hline
		$(t^3,t^8,t^{10})$ & $12$ & \makecell{$(t^3+u_1t,t^8+u_2t+u_3t^2+u_4t^4+u_5t^5+u_6t^7,$\\$t^{10}+u_7t+u_8t^2+u_9t^4+u_{10}t^5+u_{11}t^7)$}\\\hline
		$(t^4,t^5,t^{11})$ & $12$ & \makecell{$(t^4+u_1t+u_2t^2,t^5+u_3t+u_4t^2+u_5t^3+u_6t^7,$\\$t^{11}+u_7t+u_8t^2+u_9t^3+u_{10}t^6+u_{11}t^7)$}\\\hline
	\end{tabular}}
\end{table}

As in the previous section, we can apply the Restriction Method to compute the Thom Polynomials. We notice that the singularities $(t^3,t^7+t^8,0)$ and $(t^4,t^5+t^7,0)$ have a high degree term in the third coordinate of their $\mathcal{A}_e$-miniversal deformations, so if we restrict them to the $7$-jet space and $6$-jet space, respectively, that deformation direction disappears, making the principal equation computation impossible to obtain. Because of that reason, we restrict our computations to singularities of $\mathcal{A}$-codimension at most $11$, though it is possible to compute it for the quasi-homogeneous singularities of $\mathcal{A}$-codimension $12$, requiring one to obtain $102$ independent equations.

One remark that helps with the computations is that, in some cases, it is only needed to use the planar curves normal forms. For example, to compute the Thom Polynomial of codimension $7$ singularities, it is not needed to use the $(t^3,t^4,t^5)$ singularity, since the setting $$\rho_1=\alpha, \rho_2=\alpha^3\oplus\alpha^4\oplus\alpha^5$$ can also be used for the $(t^3,t^4,0)$ normal form. Furthermore, for the planar curves, there are more than one ``principal equation'', since their $\mathcal{A}_e$-miniversal deformations have a certain freedom in the powers of the toric action in the third coordinate. 

With those remarks in mind, we carry out the computations and get the following result, once again, writing the Thom Polynomials in a factorized way.

\begin{theorem}
	
The Thom Polynomials for $\mathcal{A}$-singularities of germs $(\mathbb{C},0)\to(\mathbb{C}^3,0)$ of $\mathcal{A}$-codimension less than or equal to $11$ are given in Table \ref{tab5}.

\begin{table}
	\centering
	\caption{Factorized Thom Polynomials for $\mathcal{A}\text{-cod}\leq11$ parametric space curve singularities}
	\label{tab5}
	{\renewcommand{\arraystretch}{2}\begin{tabular}{|c|c|}\hline
		\text{Normal form} & \text{Thom Polynomial}\\\hline
		$(t^2,t^3,0)$ & $c_3=-a_1^3 + a_1^2b_1 - a_1b_2 + b_3$\\\hline
		$(t^2,t^5,0)$ & $Tp_\mathcal{A}(t^2,t^3,0)(19a_1^2 - 5a_1b_1 + b_2)$\\\hline
		$(t^3,t^4,t^5)$ & $Tp_\mathcal{A}(t^2,t^3,0)(-8a_1^3 + 4a_1^2b_1 - 2a_1b_2 + b_3)
		$\\\hline
		$(t^2,t^7,0)$ & \makecell{$Tp_\mathcal{A}(t^2,t^3,0)(757a_1^4 - 336a_1^3b_1 + 35a_1^2b_1^2$\\$+ 62a_1^2b_2 - 12a_1b_1b_2 - 2a_1b_3 + b_2^2)$}\\\hline
		$(t^3,t^4,0)$  & $Tp_\mathcal{A}(t^3,t^4,t^5)(-12a_1 + b_1)$\\\hline
		$(t^3,t^5,t^7)$ & $Tp_\mathcal{A}(t^3,t^4,t^5)(37a_1^2 - 7a_1b_1 + b_2)$\\\hline
		$(t^2,t^9,0)$ & \makecell{\\$Tp_\mathcal{A}(t^2,t^3,0)(55831a_1^6 - 33113a_1^5b_1 + 6041a_1^4b_1^2$\\$+ 6309a_1^4b_2 - 343a_1^3b_1^3 - 2074a_1^3b_1b_2 - 765a_1^3b_3 + 161a_1^2b_1^2b_2$\\$+ 168a_1^2b_1b_3 + 159a_1^2b_2^2 - 7a_1b_1^2b_3 - 23a_1b_1b_2^2 - 15a_1b_2b_3$\\$+ b_1b_2b_3 + b_2^3 - b_3^2)$} \\\hline
		$(t^3,t^5,0)$ & $Tp_\mathcal{A}(t^3,t^4,t^5)(37a_1^2 - 7a_1b_1 + b_2)(-15a_1 + b_1)$\\\hline
		$(t^4,t^5,t^6)$ & $Tp_\mathcal{A}(t^3,t^4,t^5)(-27a_1^3 + 9a_1^2b_1 - 3a_1b_2 + b_3)$\\\hline
		$(t^3,t^7,t^8)$ & \makecell{$Tp_\mathcal{A}(t^3,t^4,t^5)(1489a_1^4 - 504a_1^3b_1 + 47a_1^2b_1^2 + 50a_1^2b_2$\\$- 12a_1b_1b_2 + 12a_1b_3 - b_1b_3+ b_2^2 )$}\\\hline
		$(t^4,t^5,t^7)$ & $Tp_\mathcal{A}(t^4,t^5,t^6)(-15a_1 + b_1)$\\\hline
		$(t^2,t^{11},0)$ & \makecell{\\$Tp_\mathcal{A}(t^2,t^3,0)(6159385a_1^8 - 4392432a_1^7b_1$\\$+ 1086118a_1^6b_1^2 + 830896a_1^6b_2 - 111312a_1^5b_1^3 - 372612a_1^5b_1b_2$\\$- 132142a_1^5b_3 + 4033a_1^4b_1^4 + 52160a_1^4b_1^2b_2 + 47466a_1^4b_1b_3$\\$+ 28614a_1^4b_2^2 - 2316a_1^3b_1^3b_2 - 4938a_1^3b_1^2b_3 - 7320a_1^3b_1b_2^2$\\$- 5744a_1^3b_2b_3 + 146a_1^2b_1^3b_3 + 454a_1^2b_1^2b_2^2 + 1004a_1^2b_1b_2b_3$\\$+ 304a_1^2b_2^3 + 178a_1^2b_3^2 - 38a_1b_1^2b_2b_3 - 36a_1b_1b_2^3 - 16a_1b_1b_3^2$\\$- 34a_1b_2^2b_3 + b_1^2b_3^2 + 2b_1b_2^2b_3 + b_2^4 - 2b_2b_3^2)$}\\\hline
		$(t^3,t^7+t^8,t^{11})$ & \makecell{$Tp_\mathcal{A}(t^3,t^4,t^5)(-17a_1 + b_1)(1489a_1^4 - 504a_1^3b_1 + 47a_1^2b_1^2$\\$+ 50a_1^2b_2 - 12a_1b_1b_2 + b_2^2 + 12a_1b_3 - b_1b_3)$} \\\hline
		$(t^4,t^5+t^7,t^{11})$ & $Tp_\mathcal{A}(t^4,t^5,t^6)(27a_1^2 - 2a_1b_1 + b_1^2 - 3b_2)$\\\hline
		$(t^4,t^6,t^7)$ & $Tp_\mathcal{A}(t^4,t^5,t^6)(61a_1^2 - 9a_1b_1 + b_2)$\\\hline
	\end{tabular}}
\end{table}

\end{theorem}
	
\section{Geometric method and degeneracy loci}\label{sec4}

If we treat the Thom Polynomial as a polynomial with indeterminates as the Chern classes, and the cup product as the usual product in the polynomial ring, we can factorize them. We study those factorizations and get some interesting conclusions.

\subsection{The plane case}

In Table \ref{tab3}, we can see many relations between the Thom Polynomials. One immediate remark is that every Thom Polynomial is divisible by $Tp_\mathcal{A}(A_2)$. Furthermore, the Thom Polynomial for $E_{k}$ and $W_{k}$ are also divisible by $Tp_\mathcal{A}(E_6)$, and the Thom Polynomial for $W_k$ singularities are divisible by $Tp_\mathcal{A}(W_{12}')$.

This could be explained by the projection formula, since if we consider the singular locus $\Sigma_{A_2}$ and take a worse singularity $\eta$, then its singular locus $\Sigma_\eta$ is contained in $\Sigma_{A_2}$, so we consider the inclusion $\iota:\Sigma_{A_2}\to J(1,2)$, and denoting by $\alpha$ the cohomology class representing the $\eta$ singularity, we get $$\iota_\ast(\alpha)=\iota_\ast(1\cdot\iota^\ast(\beta))=\iota_\ast(1)\cdot\beta,$$ and $\iota_\ast(1)$ is $Tp_\mathcal{A}(A_2)$, provided that $\alpha=\iota^\ast(\beta)$ is in the image of $\iota^\ast$, and to prove this, probably it is needed to use the smoothness of $\Sigma_{A_2}$, since we do not have divisibility between $Tp_\mathcal{A}(A_6)$ and $Tp_\mathcal{A}(A_8)$, for example. The locus of $A_2$ singularity is simply the vanishing of the $1$-jet, so it is smooth, and the same explanation applies to $E_6$ (the vanishing of the $2$-jet) and to $W_{12}'$ (the vanishing of the $3$-jet). 

Another explanation could be the Thom-Porteous Formula, which we write it here and we refer the reader to see \cite{5} for more details.

\begin{theorem}[Thom-Porteous Formula]
	
For a suitable generic vector bundle map $f:E\to F$, where $E$ and $F$ are vector bundles over a smooth manifold $M$ of rank $m$ and $n$, respectively, then the degeneracy loci $$D^k(f)=\{x\in M\mid \textnormal{dim}(\textnormal{ker}(f_x))\geq k\}$$ is expressed by $$\iota_\ast([D^k(f)])=\textnormal{det}([c_{k+\kappa-i+j}(F-E)]_{1\leq i,j\leq k})\frown[M],$$ where $\kappa=n-m$ is the relative codimension of $f$.
	
\end{theorem}

The Thom-Porteous Formula can be used to determine some of such factors. For example, let us describe how to use this formula to obtain the Thom Polynomial for the $\mathcal{A}$-singularity $A_4=(t^2,t^5)$. First, we write an element in the jet space as \begin{equation}\label{exp1}\sigma=\mathbb{a}_1t+\mathbb{a}_2t^2+\cdots=(a_{11}t+a_{12}t^2+a_{13}t^3+\cdots,a_{21}t+a_{22}t^2+a_{23}t^3+\cdots),\end{equation} where $\mathbb{a}_i=(a_{1i},a_{2i})$. Note that there is no linear terms in the $A_4$ normal form, so we will consider our vector bundles over the singular locus of type $A_2$, $\Sigma_{A_2}\subset J(1,2)$, in order to have $\mathbb{a}_1=0$.

Now, since the $A_4$ normal form starts with a $t^2$ term, we must have that $\mathbb{a}_2\neq0$. Following the technique shown in \cite{7} by Arnol'd, through a change of coordinates in the target, we can eliminate the term $t^2$ from one of the two components of $\sigma$. For example, if $a_{12}\neq0$, then we can consider $$k(x, y)=\left(x, y-\frac{a_{22}}{a_{12}}x\right).$$ We have that $k$ is a change of coordinates in the target such that, when composed with $\sigma$, the term $t^2$ in the second coordinate vanishes. If $a_{22}\neq0$, the argument is similar.

In order to have no $t^3$ terms, it is needed to eliminate $\mathbb{a}_3$ with $\mathbb{a}_2$, so we require that $\mathbb{a_3}$ is parallel to $\mathbb{a}_2$. Indeed, suppose that we have $\mathbb{a}_3=\lambda\mathbb{a}_2$, with $\lambda\in\mathbb{C}$. For any $\alpha\in\mathbb{C}$, consider the change of coordinates in the source $$h^{-1}(t)=t-\alpha t^2.$$ Then, composing $\sigma$ with $h^{-1}$, we get $$(\sigma\circ h^{-1})(t)=(a_{12}t^2+(-2a_{12}\alpha+a_{13})t^3+o(t^4), a_{22}t^2+(-2a_{22}\alpha+a_{23})t^3+o(t^4)),$$ where $o(t^4)$ denotes terms $t^k$, with $k\geq4$. Since $a_{i3}=\lambda a_{i2}$, for $i=1,2$, then the coefficients of $t^3$ in $(\sigma\circ h^{-1})(t)$ are $$-2a_{i2}\alpha+\lambda a_{i2}=(-2\alpha+\lambda)a_{i2},\text{ for }i=1,2.$$ Therefore, if we take $\alpha=\lambda/2$, we eliminate the term $t^3$.

With those requirements, the jet $\sigma$ is equivalent to one that starts with $(t^2,t^5)$, since now we can eliminate the term $t^4$ in both coordinates by using $t^2$ through $$h^{-1}(t)=t-\beta t^3 \text{ and }k(x,y)=(x, y-\gamma x^2),$$ for some appropriate $\beta$ and $\gamma$. In the end, we finally get $A_4$, since even power terms $t^{2k}$ can be eliminated by using $t^2$ and odd power terms $t^{2k+1}$ higher than $t^5$ can be eliminated by using both $t^2$ and $t^5$ in a similar way as above.

Geometrically, let us consider two vector bundles, $E$ and $F$, of rank $1$ and $2$, respectively, over $\Sigma_{A_2}$, where the former will parameterize the powers $t^i$ and the latter will parameterize the coefficient vectors $\mathbb{a}_i$. With the above discussion, we consider the mapping $$f:E^{\otimes2}\oplus E^{\otimes3}\to F$$ with $k=1$, which encodes the condition that $\{\mathbb{a}_2,\mathbb{a}_3\}$ is linearly dependent. Then, we get from the Thom-Porteous Formula that $$[D^k(f)]=\left|\begin{array}{c}c_2(F-E)\end{array}\right|=-5a_1+b_1$$ and we have that $$Tp_\mathcal{A}(A_4)=(a_1^2-a_1b_1+b_2)(-5a_1+b_1).$$

In order to apply this explanation to prove the divisibility property in general, we face the problem to determine the suitable vector bundle map and the degeneracy loci corresponding to the singularity in order to compute its Thom Polynomial, though it is already known that we can use the Thom-Porteous formula to compute the Thom Polynomial of certain singularities by considering vector bundles over specific loci that are smooth and one mapping between them.

This divisibility property can be useful if one can find a theoretical way to compute the expression for $\beta$. Another usefulness is for concrete computations, since $Tp_\mathcal{A}(A_2)$ is already determined, so to compute $Tp_\mathcal{A}(A_6)$, for example, we can write it as $$Tp_\mathcal{A}(A_6)=(a_1^2-a_1b_1+b_2)(x_1a_1^2+x_2a_1b_1+x_3b_1^2+x_4b_2),$$ and apply the Restriction Method to determine those four unknown integers, in contrast to nine when using the usual Restriction Method.

There are some generalizations of Chern classes formulas for degeneracy loci, such as the works made originally by Kempf and Laksov in \cite{19}, and then expanded further by Fehér and Rimányi in \cite{17} and by Anderson and Fulton in \cite{9}. Inspired by the application of the Thom-Porteous Formula, we managed to find a connection between the formulas shown in \cite{9} and the Thom Polynomials that we computed. The advantage of this formula is that we can consider a finite number of mappings, so we can encode many conditions at once, in contrast to the Thom-Porteous Formula, where many conditions must be encoded in one single mapping, which can be impossible to achieve for worse singularities, as their geometry gets more complicated.

To explain how to apply the formula obtained by Anderson and Fulton, let us first state their framework. Consider a sequence of vector bundles 

\begin{tikzcd}E_{p_1}\arrow[hook]{r}& E_{p_2}\arrow[hook]{r}& \cdots\arrow[hook]{r}& E_{p_s}\arrow{r}{\varphi}& F_{q_1}\arrow[two heads]{r}& F_{q_2}\arrow[two heads]{r}& \cdots\arrow[two heads]{r}& F_{q_s}\end{tikzcd} on a variety $X$, where the subscripts indicate the ranks. Then, for each $i$, there is an induced mapping $h_i:E_{p_i}\to F_{q_i}$. Let $\textnormal{\textbf{p}}=(p_1,\ldots,p_s)$ and $\textnormal{\textbf{q}}=(q_1,\ldots,q_s)$ be the rank vectors. Given a vector $\textnormal{\textbf{r}}=(r_1,\ldots,r_s)$ of non-negative integers, put $\tau=(\textnormal{\textbf{r}}, \textnormal{\textbf{p}}, \textnormal{\textbf{q}})$ and define the degeneracy locus corresponding to $\tau$ as $$\Omega_\tau=\{x\in X\mid \text{rank}(E_{p_i}(x)\to F_{q_i}(x))\leq r_i, \text{ for }1\leq i\leq s\}.$$ Put $k_i=p_i-r_i$, $l_i=q_i-r_i$ and suppose that \begin{align*} 0 < p_1\leq p_2 \leq \cdots \leq p_s,\\
q_1 \geq q_2 \geq \cdots \geq q_s>0,\\
0 < k_1 < k_2 < \cdots < k_s,\\
l_1 \geq l_2 \geq \cdots \geq l_s \geq 0.\end{align*} Define $\lambda_{k_i}=l_i$, and whenever there is a gap in the sequence $(k_i)_{1 \leq i \leq s}$, fill the remaining parts in minimally such that $$\lambda_1 \geq \lambda_2 \geq \cdots \geq \lambda_{k_s} \geq 0.$$ Finally, putting $c(k_i)=c(F_{q_i}-E_{p_i})$, completing the sequence with $c(k)=c(k_i)$, for $k_{i-1} < k \leq k_i$, and denoting by $\ell=k_s$, Anderson and Fulton proved that the fundamental class associated to $\Omega_\tau$ is expressed by the $\ell\times\ell$ determinant $$[\Omega_\tau]=\left|\begin{array}{cccc}c(1)_{\lambda_1} & c(1)_{\lambda_1+1} & \cdots & c(1)_{\lambda_1+\ell-1}\\c(2)_{\lambda_2-1} & c(2)_{\lambda_2} & \cdots & c(2)_{\lambda_2+\ell-2}\\\vdots & \vdots & \ddots & \vdots\\c(\ell)_{\lambda_\ell+1-\ell} & c(\ell)_{\lambda_\ell+2-\ell} & \cdots & c(\ell)_{\lambda_\ell}\end{array}\right|,$$ where $c(i)_j$ is the $j$-th quotient Chern class that appears in $c(i)$.

Let us apply this formula to an example. Consider the $E_8$ singularity in the plane case, with normal form $(t^3,t^5)$. Since we can encode many conditions simultaneously, we can consider our vector bundles $E$ and $F$ over some complex projective space $\mathbb{P}_\mathbb{C}^N$, with ranks $1$ and $2$, respectively.

Similar to what we did to obtain $Tp_\mathcal{A}(A_4)$ using the Thom-Porteous Formula, given a jet $\sigma$ as in \ref{exp1}, it is needed that there is neither $t$ nor $t^2$, so we require that $\mathbb{a}_1=\mathbb{a}_2=0$ and, thus, it starts with the term $t^3$. Now, we can eliminate $a_{23}$ with $a_{13}$ to eliminate the term $t^3$ in the second component, and we require that $\{\mathbb{a}_3,\mathbb{a}_4\}$ is linearly dependent in order to have no $t^4$ terms. With those requirements, the jet $\sigma$ is equivalent to one that starts with $(t^3,t^5)$, and its higher order terms can be eliminated to get $E_8$.

With the above requirements, we consider a sequence of fiber bundles and morphisms

$$\begin{tikzcd}
	E\oplus E^{\otimes2}\arrow[hook]{r} & E\oplus E^{\otimes 2}\oplus E^{\otimes 3}\oplus E^{\otimes 4}\arrow{r} & F \arrow[two heads]{r}& F.
\end{tikzcd}$$

We have that the respective rank vectors are $\textnormal{\textbf{p}}=(p_1,p_2)=(2,4)$ and $\textnormal{\textbf{q}}=(q_1,q_2)=(2,2)$. Now we consider the morphisms $h_i:E_{p_i}\to F_{q_i}$ that translate all of the requirements that we put on the coefficients $\mathbb{a}_i$ as follows: \begin{align*}
h_1:E\oplus E^{\otimes 2}\to F,\text{ }r_1=0,\text{ }k_1=2,\text{ }l_1=2,\\
h_2:E\oplus E^{\otimes 2}\oplus E^{\otimes 3}\oplus E^{\otimes 4}\to F,\text{ }r_2=1,\text{ }k_2=3,\text{ }l_2=1,
\end{align*} where $r_i$ is the rank of the morphism $h_i$. Then, as described above, we denote by $\lambda_2=l_1=2$ and $\lambda_3=l_2=1$. Since there is a gap in the sequence of the kernel dimensions $(k_1,k_2)=(2,3)$, we follow Fulton's method and fill it in minimally with $\lambda_1=l_1=2$ and $c(1)=c(2)$, where $$c(2)=c(F-E\oplus E^{\otimes2})\textnormal{ and }c(3)=c(F-E\oplus E^{\otimes2}\oplus E^{\otimes3}\oplus E^{\otimes4})$$ are the total quotient Chern classes of each corresponding mapping $h_1$ and $h_2$. Then, for $\tau=(\textnormal{\textbf{r}},\textnormal{\textbf{p}},\textnormal{\textbf{q}})$, where $\textnormal{\textbf{r}}=(r_1,r_2)=(0,1)$, the degeneracy locus $$\Omega_\tau=\{x\in \mathbb{P}_\mathbb{C}^N\mid \textnormal{rank}(E_{p_i}(x)\to F_{q_i}(x))\leq r_i,\textnormal{ for }i=1,2\}$$ has cohomology class \begin{align*}
\left[\Omega_\tau\right]&=\left|\begin{array}{ccc}c(1)_2&c(1)_3&c(1)_4\\c(2)_1&c(2)_2&c(2)_3\\0&1&c(3)_1\end{array}\right|\\&=\makecell{-28a_1^5 + 46a_1^4b_1 - 20a_1^3b_1^2 - 35a_1^3b_2 + 2a_1^2b_1^3\\+ 26a_1^2b_1b_2- 3a_1b_1^2b_2- 7a_1b_2^2 + b_1b_2^2}\\&=Tp_\mathcal{A}(E_8).
\end{align*}

Another advantage of this method is the freedom in the choices that the settings can be defined, provided that it is consistent with the conditions given in \cite{9}. For example, we can take a third mapping to encode the condition $\mathbb{a}_1=0$ separately, that is, consider \begin{align*}h_1:E\to F,\text{ }r_1=0,\text{ }k_1=1,\text{ }l_1=2.\\
h_2:E\oplus E^{\otimes 2}\to F,\text{ }r_2=0,\text{ }k_2=2,\text{ }l_2=2,\\
h_3:E\oplus E^{\otimes 2}\oplus E^{\otimes 3}\oplus E^{\otimes 4}\to F,\text{ }r_3=1,\text{ }k_3=3,\text{ }l_3=1,\end{align*}
With this configuration, there is no gap in the sequence $(k_1,k_2,k_3)=(1,2,3)$, but we get the same expression that coincides with $Tp_\mathcal{A}(E_8)$. 

Also, if we consider the vector bundles over the $A_2$ locus instead of the complex projective space, we consider only the conditions $\mathbb{a}_2=0$ and $\{\mathbb{a}_3,\mathbb{a}_4\}$ spanning an $1$-dimensional vector space, that is, \begin{align*}
h_1:E^{\otimes2}\to F,\text{ }r_1=0,\text{ }k_1=1,\text{ }l_1=2,\\
h_2:E^{\otimes 2}\oplus E^{\otimes3}\oplus E^{\otimes4}\to F,\text{ }r_2=1,\text{ }k_2=2,\text{ }l_2=1.
\end{align*} Then there is also no gap in the sequence $(k_1, k_2) = (1, 2)$ and we obtain $$\left[\Omega_\lambda(\tau)\right]=\left|\begin{array}{cc} c(1)_2 & c(1)_3\\
1 & c(2)_1\end{array}\right|=(4a_1^2-2a_1b_1+b_2)(-7a_1+b_1)$$ and $$Tp_\mathcal{A}(A_2)\cdot\left[\Omega_\lambda(\tau)\right]=Tp_\mathcal{A}(E_8).$$ Similarly, we can consider the vector bundles over the $E_6$ locus, and the degeneracy locus expression will be the factor $-7a_1+b_1$, that when multiplied by $Tp_\mathcal{A}(E_6)$ results in $Tp_\mathcal{A}(E_8)$.

In Table \ref{tab6}, we list how to apply the above method for some singularities and one possible configuration for each of them, considering the vector bundles over some complex projective space, though one can change the settings like we did in the above example, and it is also possible to simplify them by considering the vector bundles over a smooth locus, which will result in a degeneracy loci expression that must be multiplied by the corresponding Thom Polynomial to get the desired expression.

\begin{table}
	\centering
	\caption{Some configurations for applying the degeneracy loci formula}
	\label{tab6}
	{\renewcommand{\arraystretch}{2}\begin{tabular}{|c|c|}\hline
	\textnormal{Symbol} & \textnormal{Configuration}\\\hline
	$A_2$ & $h_1:E\to F,\text{ }r_1=0,\text{ }k_1=1, \text{ }l_1=2$\\\hline
	$A_4$ & \makecell{$h_1:E\to F,\text{ }r_1=0,\text{ }k_1=1,\text{ }l_1=2$\\$h_2:E\oplus E^{\otimes2}\oplus E^{\otimes3}\to F,\text{ }r_2=1,\text{ }k_2=2,\text{ }l_2=1$}\\\hline
	$E_6$ & $h_1:E\oplus E^{\otimes2}\to F,\text{ }r_1=0,\text{ }k_1=2,\text{ }l_1=2$\\\hline
	$E_8$ & \makecell{$h_1:E\oplus E^{\otimes 2}\to F,\text{ }r_1=0,\text{ }k_1=2,\text{ }l_1=2$\\$h_2:E\oplus E^{\otimes 2}\oplus E^{\otimes 3}\oplus E^{\otimes 4}\to F,\text{ }r_2=1,\text{ }k_2=3,\text{ }l_2=1$}\\\hline
	$E_{12}'$ & \makecell{$h_1:E\oplus E^{\otimes2}\to F,\text{ }r_1=0,\text{ }k_1=2,\text{ }l_1=2$\\$h_2:E\oplus E^{\otimes2}\oplus E^{\otimes3}\oplus E^{\otimes4}\oplus E^{\otimes5}\to F,\text{ }r_2=1,\text{ }k_2=4,\text{ }l_2=1$}\\\hline
	$W_{12}'$ & $h_1:E\oplus E^{\otimes2}\oplus E^{\otimes3}\to F,\text{ }r_1=0,\text{ }k_1=3,\text{ }l_1=2$\\\hline
	$W_{12}$ & \makecell{$h_1:E\oplus E^{\otimes2}\oplus E^{\otimes3}\to F,\text{ }r_1=0,\text{ }k_1=3,\text{ }l_1=2$\\$h_2:E\oplus E^{\otimes2}\oplus E^{\otimes3}\oplus E^{\otimes4}\oplus E^{\otimes7}\to F,\text{ }r_2=1,\text{ }k_2=4,\text{ }l_2=1$}\\\hline
\end{tabular}}
\end{table}

Unfortunately, the author could not proceed further to get $Tp_\mathcal{A}(A_6)$ and beyond directly from the degeneracy loci formula. Naturally, it could be thought that $A_6$ is determined by $\mathbb{a}_1=0$ and $\{\mathbb{a}_2,\mathbb{a}_3,\mathbb{a}_5\}$ spanning a $1$-dimensional vector space, so we could eliminate $\mathbb{a}_3$ and $\mathbb{a}_5$ with $\mathbb{a}_2$. Now, considering the vector bundles over the complex projective space and the setting \begin{align*}
	h_1:E\to F,\text{ }r_1=0,\text{ }k_1=1,\text{ }l_1=2,\\
	h_2:E\oplus E^{\otimes2}\oplus E^{\otimes3}\oplus E^{\otimes5}\to F,\text{ }r_2=1,\text{ }k_2=3,\text{ }l_2=1,
\end{align*} we get the equality $$[\Omega_\tau]=Tp_\mathcal{A}(A_6)-Tp_\mathcal{A}(E_6),$$ and other settings also result in the same equality. The meaning of this equality is that the $\mathcal{A}$-codimension of the $A_6$ and the $E_6$ singularities are the same and the support of that degeneracy locus is the union of those two loci. Indeed, the mapping $h_1$ says that $\mathbb{a}_1=0$, so the jet $\sigma$ starts with $t^2$ or higher order terms. The mapping $h_2$ says that it is possible to eliminate $t^5$ by using either $t^2$ or $t^3$. However, it does not control if $\mathbb{a}_2\neq0$ or if $\mathbb{a}_3\neq0$. If $\mathbb{a}_2\neq0$, then we get $A_6$, and if $\mathbb{a}_2=0$ and $\mathbb{a}_3\neq0$, then we get $E_6$. Furthermore, it is not possible to consider a separate mapping  to encode only the condition $\mathbb{a}_2\neq0$. For instance, if we append the mapping $$h: E\oplus E^{\otimes2}\to F,\text { }r=1, \text{ }k=1, \text{ }l=1$$ to the above setting, the vector of kernel dimensions will be $(1, 1, 3)$, which is not a strictly increasing sequence of integers, so it is not possible to apply the degeneracy loci formula.

Recently, Fukui and Hoshino made a complete $\mathcal{A}$-recognition of parametric curves germs in the plane in terms of the intrinsic derivatives in \cite{23}. Therefore, their results and techniques can be used to obtain further Thom Polynomials if this application of the degeneracy loci formula is well understood.

\subsection{The space case}

As we did in the plane curves case, we can apply the formula in \cite{9} to obtain some Thom Polynomials. In this case, we consider $F$ a vector bundle of rank $3$ over some complex projective space. In Table \ref{tab7}, we write one configuration for some of the singularities that we have considered.

\begin{table}
	\centering
	\caption{Some configurations for applying the degeneracy loci formula}
	\label{tab7}
	{\renewcommand{\arraystretch}{2}\begin{tabular}{|c|c|}\hline
		\textnormal{Normal Form} & \textnormal{Configuration}\\\hline
		$(t^2,t^3,0)$ & $h_1:E\to F,\text{ }r_1=0,\text{ }k_1=1, \text{}l_1=3$\\\hline
		$(t^2,t^5,0)$ & \makecell{$h_1:E\to F,\text{ }r_1=0,\text{ }k_1=1, \text{}l_1=3$\\$h_2:E\oplus E^{\otimes2}\oplus E^{\otimes3}\to F,\text{ }r_2=1,\text{ }k_2=2, \text{}l_2=2$}\\\hline
		$(t^3,t^4,t^5)$ & $h_1:E\oplus E^{\otimes2}\to F,\text{ }r_1=0,\text{ }k_1=2\text{ }l_1=3$\\\hline
		$(t^3,t^4,0)$ & \makecell{$h_1:E\oplus E^{\otimes2}\to F,\text{ }r_1=0,\text{ }k_1=2\text{ }l_1=3$\\$h_2:E\oplus E^{\otimes2}\oplus E^{\otimes3}\oplus E^{\otimes4}\oplus E^{\otimes5}\to F,\text{ }r_2=2,\text{ }k_2=3,\text{ }l_2=1$}\\\hline
		$(t^3,t^5,t^7)$ & \makecell{$h_1:E\oplus E^{\otimes2}\to F,\text{ }r_1=0,\text{} k_1=2,\text{ }l_1=3$\\$h_2:E\oplus E^{\otimes2}\oplus E^{\otimes3}\oplus E^{\otimes4}\to F,\text{ }r_2=1,\text{ }k_2=3,\text{ }l_2=2$}\\\hline
		$(t^3,t^7,t^8)$ & \makecell{$h_1:E\oplus E^{\otimes2}\to F,\text{ }r_1=0,\text{} k_1=2,\text{ }l_1=3$\\$h_2:E\oplus E^{\otimes2}\oplus E^{\otimes3}\oplus E^{\otimes4}\oplus E^{\otimes5}\to F,\text{ } r_2=1,\text{ }k_2=4,\text{ }l_2=2$}\\\hline
		$(t^4,t^5,t^6)$ & $h_1:E\oplus E^{\otimes2}\oplus E^{\otimes3}\to F,\text{ }r_1=0,\text{ }k_1=3,\text{ }l_1=3$\\\hline
		$(t^4,t^5,t^7)$ & \makecell{$h_1:E\oplus E^{\otimes2}\oplus E^{\otimes3}\to F,\text{ }r_1=0,\text{ }k_1=3,\text{ }l_1=3$\\$h_2:E\oplus E^{\otimes2}\oplus E^{\otimes3}\oplus E^{\otimes4}\oplus^{\otimes5}\oplus E^{\otimes6}\to F,\text{ }r_2=2,\text{ }k_2=4,\text{ }l_2=1$}\\\hline
		$(t^4,t^6,t^7)$ & \makecell{$h_1:E\oplus E^{\otimes2}\oplus E^{\otimes3}\to F,\text{ }r_1=0,\text{ }k_1=3,\text{ }l_1=3$\\$h_2:E\oplus E^{\otimes2}\oplus E^{\otimes3}\oplus E^{\otimes4}\oplus E^{\otimes5}\to F,\text{ }r_2=1,\text{ }k_2=4,\text{ }l_2=2$}\\\hline
	\end{tabular}}
\end{table}

Furthermore, similarly to the plane case, if we consider $[\Omega_\tau]$ the degeneracy loci expression coming from $\mathbb{a}_1=0$ and $\{\mathbb{a}_2,\mathbb{a}_3,\mathbb{a}_5\}$ spanning a $1$-dimensional vector space, we also have the equality $$[\Omega_\tau]=Tp_\mathcal{A}(t^2,t^7,0)-Tp_\mathcal{A}(t^3,t^4,0).$$ In the space case, we have some new singularities, for instance, $(t^3,t^4,t^5)$. The singularity $(t^3,t^4,0)$ is determined from $(t^3,t^4,t^5)$ by requiring that $\mathbb{a}_5$ is a linear combination of $\mathbb{a}_3$ and $\mathbb{a}_4$, so we consider the degeneracy loci coming from $\mathbb{a}_1=\mathbb{a}_2=0$ and $\{\mathbb{a}_3,\mathbb{a}_4,\mathbb{a}_5\}$ spanning a $2$-dimensional vector space, obtaining $$[\Omega_\tau]=Tp_\mathcal{A}(t^3,t^4,0).$$ Also, if we consider the degeneracy loci coming from $\mathbb{a}_1=\mathbb{a}_2=0$, $\{\mathbb{a}_3,\mathbb{a}_4\}$ spanning a $1$-dimensional vector space and $\{\mathbb{a}_3,\mathbb{a}_5,\mathbb{a}_7\}$ spanning a $2$-dimensional vector space, we get the equality $$[\Omega_\tau]=Tp_\mathcal{A}(t^3,t^5,0)-Tp_\mathcal{A}(t^4,t^5,t^6).$$

\section{Thom Polynomial for some parametric curves germs in $\mathbb{C}^N$}

We saw in the previous sections that, for explicit dimensions and singularities, the Thom Polynomial can be computed with enough time and computational power. In \cite{7}, Arnol'd obtains the classification of simple singularities of parametric curves in $\mathbb{C}^N$. In this case, the computation of the Thom Polynomial can no longer be made by the Restriction Method, or at least, not explicitly enough to be possible to use a mathematical software.

For example, following the notation in \cite{7}, consider the singularity $(2,3)^N$, whose upper index denotes its $\mathcal{A}$-codimension, that is, the curve $(t^2,t^3,0,\ldots,0)$ in $\mathbb{C}^N$. In the previous sections, we obtained its Thom Polynomial for $N=2,3$, and we could carry out its computations for any concrete $N$. Now, for an arbitrary $N$, it is still possible to compute it with the restriction method, which is done by Rimányi in \cite{3} for the $\mathcal{K}$-equivalence case, so we just substitute the quotient Chern classes by the expanded expressions in Chern classes of the source and those of the target and we get $$Tp_\mathcal{A}((t^2,t^3)^N)=\sum_{i=0}^N(-1)^{N-i}a_1^{N-i}b_i.$$ Though this expression is already known due to Thom (\cite{13}) and Porteous (\cite{12}). Furthermore, the Thom Polynomial for the $\mathcal{K}$-singularity $A_2$, that corresponds to $(3,4,5)^{2N}$ in our case, is also already known due to Ronga (\cite{14}).

The above expression shows that the task proposed in this section is already difficult enough for the most basic singularity. Even for concrete examples, the required computational power quickly becomes enormous. For example, the Thom Polynomial for $(3,4,5)^8=(t^3,t^4,t^5,0)$ requires us to obtain and solve a linear system consisting of $94$ equations.

For $\mathcal{K}$-equivalence, in \cite{11}, Feh\'{e}r and Rim\'{a}nyi obtained a series expression, called the Thom Series, and in \cite{8}, Kazarian obtained an interesting formula for the general Thom Polynomial using non-associative Hilbert schemes to find a suitable resolution of singularities. It is a well-known fact that if the $\mathcal{A}$-orbit of an $\mathcal{A}$-singularity is dense in its $\mathcal{K}$-orbit, then its Thom Polynomial can be written in quotient Chern classes. Therefore, using Kazarian's result and substituting the quotient Chern classes, it is possible to recover the above formulas, though it only works for very few $\mathcal{A}$-singularities, namely $(2,3)^N$, $(3,4,5)^{2N}$ and $(4,5,6,7)^{3N}$, because of the denseness condition of the orbits.

One problem is to generalize Kazarian's method to $\mathcal{A}$-equivalence without requiring the denseness condition of the $\mathcal{A}$-orbit. Clearly, the associated algebra must be replaced by something that can detect higher powers. For example, $(t^2,t^3)^N$ and $(t^2,t^5)^{2N-1}$ have the same local algebra, so Kazarian's method will give the same generating series and the same Thom Polynomial for both singularities, which is not true.

We now propose, following the method using the degeneracy loci Chern classes formula. For simplicity, we denote the total quotient Chern class $$c(i_1,i_2,\ldots,i_k)=c(F-E^{\otimes i_1}\oplus E^{\otimes i_2}\oplus\cdots\oplus E^{\otimes i_k}),$$ where $F$ is a vector bundle of rank $N$ and $E$ is a vector bundle of rank $1$ over some complex projective space, and we denote by $c(i_1,i_2,\ldots,i_k)_j$ the $j$-th quotient Chern class that appears in the expression of the total quotient Chern class $c(i_1,i_2,\ldots,i_k)$.

\begin{conjecture}

The Thom Polynomials for some singularities of parametric curves in $\mathbb{C}^N$ are presented in Table \ref{tab8}. 

\begin{table}
	\centering
	\caption{Proposed Thom Polynomial for parametric curve singularities in $\mathbb{C}^N$}
	\label{tab8}
	{\renewcommand{\arraystretch}{2}\begin{tabular}{|c|c|}\hline
		\textnormal{Normal form} & \textnormal{Thom Polynomial}\\\hline
		$(2,3)^N$ & $\sum_{i=0}^N(-1)^{N-i}a_1^{N-i}b_i=c_N$\\\hline
		$(3,4,5)^{2N}$ & \makecell{$c_N^2+\sum_{k=1}^N2^{k-1}c_N\left(\sum_{i=0}^{N-k}(-1)^{N-i}a_1^{N-i}b_i\right)$\\$=c(1,2)_N^2-c(1,2)_{N-1}c(1,2)_{N+1}$}\\\hline
		$(4,5,6,7)^{3N}$ & $\left|\begin{array}{ccc}c(1,2,3)_N & c(1,2,3)_{N+1} & c(1,2,3)_{N+2}\\c(1,2,3)_{N-1} & c(1,2,3)_N & c(1,2,3)_{N+1}\\c(1,2,3)_{N-2} & c(1,2,3)_{N-1} & c(1,2,3)_N\end{array}\right|$\\\hline
		$(2,5)^{2N-1}$ & $c(1)_Nc(1,2,3)_{N-1}-c(1)_{N+1}c(1,2,3)_{N-2}$\\\hline
		$(3,4)^{3N-2}$ & $\left|\begin{array}{ccc}c(1,2)_N & c(1,2)_{N+1} & c(1,2)_{N+2}\\c(1,2)_{N-1} & c(1,2)_N & c(1,2)_{N+1}\\c(1,2,3,4,5)_{N-4} & c(1,2,3,4,5)_{N-3} & c(1,2,3,4,5)_{N-2}\end{array}\right|$\\\hline
		$(3,5,7)^{3N-1}$ & $\left|\begin{array}{ccc}c(1,2)_N & c(1,2)_{N+1} & c(1,2)_{N+2}\\c(1,2)_{N-1} & c(1,2)_N & c(1,2)_{N+1}\\c(1,2,3,4)_{N-3} & c(1,2,3,4)_{N-2} & c(1,2,3,4)_{N-1}\end{array}\right|$\\\hline
		$(4,5,6)^{4N-3}$ & \makecell{$\left|\begin{array}{cccc} c(1,2,3)_N & c(1,2,3)_{N+1} & c(1,2,3)_{N+2}  & c(1,2,3)_{N+3}  \\ c(1,2,3)_{N-1} & c(1,2,3)_N & c(1,2,3)_{N+1} & c(1,2,3)_{N+2} \\ c(1,2,3)_{N-2} & c(1,2,3)_{N-1} & c(1,2,3)_N & c(1,2,3)_{N+1} \\ d_{N-6} & d_{N-5} & d_{N-4} & d_{N-3}\end{array}\right|$\\\text{where }$d=c(1,2,3,4,5,6,7)$}\\\hline
	\end{tabular}}
\end{table}
	
\end{conjecture}	
	
If, as seen in the previous sections, $Tp_\mathcal{A}((2,7)^{3N-2})-Tp_\mathcal{A}((3,4)^{3N-2})$ is equal to $$\left|\begin{array}{ccc}c(1)_N & c(1)_{N+1} & c(1)_{N+2}\\ c(1,2,3,5)_{N-2} & c(1,2,3,5)_{N-1} & c(1,2,3,5)_N\\ c(1,2,3,5)_{N-3} & c(1,2,3,5)_{N-2} & c(1,2,3,5)_{N-1} \end{array}\right|$$ and $Tp_\mathcal{A}((3,5)^{4N-3})-Tp_\mathcal{A}((4,5,6)^{4N-3})$ is equal to $$\left|\begin{array}{cccc}c(1,2)_N & c(1,2)_{N+1} & c(1,2)_{N+2} & c(1,2)_{N+3}\\c(1,2)_{N-1} & c(1,2)_N & c(1,2)_{N+1} & c(1,2)_{N+2}\\c(1,2,3,4)_{N-3} & c(1,2,3,4)_{N-2} & c(1,2,3,4)_{N-1} & c(1,2,3,4)_N\\c(1,2,3,4,5,7)_{N-5} & c(1,2,3,4,5,7)_{N-4} & c(1,2,3,4,5,7)_{N-3} & c(1,2,3,4,5,7)_{N-2}\end{array}\right|,$$ for every $N$, then we can go one step further and conjecture the general Thom Polynomial expression for $(2,7)^{3N-2}$ and $(3,5)^{4N-3}$ as the sum of two determinants.

\backmatter

\bmhead{Acknowledgments}

The author wishes to extend his heartfelt appreciation to FAPESP (S\~{a}o Paulo Research Foundation) for their invaluable partial support, process numbers 2020/14442-9 and 2023/03086-5, which played an important role in the development of this scientific research, and to Toru Ohmoto, for many fruitful discussions and key remarks that allowed the discovery of the results written in this paper while the author was in his abroad internship granted by FAPESP to study at Waseda University.

\bibliography{sn-bibliography}

\end{document}